\documentclass{article}

\usepackage{amsfonts,amssymb,amsmath}
\usepackage[dvips]{graphicx}
\usepackage[pdftex]{hyperref}
\usepackage{afterpage}
\usepackage{algorithm}
\usepackage{algorithmic}
\usepackage{color}
\usepackage{fancyhdr}
\usepackage{natbib}

\newcommand{\mean}{\mathbb{E}}
\newcommand{\nelite}{\ensuremath N^{\rm elite}}
\newcommand{\estimatorvalue}{\hat{\phi}}
\newcommand{\estimatorrv}{\widehat{\Phi}}
\newcommand{\targetrv}{\Phi}

\newcommand{\var}{\mathbb{V}}
\newcommand{\xith}[2]{\ensuremath #1^{(#2)}}

\begin{document}
\title{Bias-Variance Techniques for Monte Carlo Optimization:\\Cross-validation for the CE Method }
\author{Dev Rajnarayan and David Wolpert}
\date{\today}
\maketitle

\section{Introduction}
\label{sec:intro}
In this paper, we examine the CE method in the broad context of Monte
Carlo Optimization (MCO)~\citep{erno98, roca04} and Parametric
Learning (PL), a type of machine learning.  A well-known overarching
principle used to improve the performance of many PL algorithms is the
bias-variance tradeoff~\citep{wolp97}. This tradeoff has been used to
improve PL algorithms ranging from Monte Carlo estimation of
integrals~\citep{lepa78}, to linear estimation, to general statistical
estimation~\citep{brei96,brei96a}.  Moreover, as described by
\citet{wora07}, MCO is  very closely related to PL. Owing to this similarity, 
 the bias-variance tradeoff affects MCO performance, just as it does 
 PL performance.

In this article, we exploit the bias-variance
tradeoff to enhance the performance of MCO algorithms. We use
the technique of cross-validation, a technique based on the bias-variance tradeoff, to
significantly improve the performance of the Cross Entropy (CE)
method, which is an MCO algorithm.  In previous work we have confirmed
that other PL techniques improve the perfomance of other MCO
algorithms \cite[see][]{wora07}. We conclude that the many techniques
pioneered in PL could be investigated as ways to improve MCO
algorithms in general, and the CE method in particular.

The rest of the paper is organized as follows. In
Sec.~\ref{sec:BPVDerivation}, we present an overview of the
bias-variance tradeoff. In Sec.~\ref{sec:applications}, we describe a
few ways to exploit this tradeoff, starting from the relatively simple
case of Monte Carlo integration, and proceeding to the more complex
case of MCO. We also describe the original exploitation of this
tradeoff, as a way to improve PL algorithms.  In
Sec.~\ref{sec:PLMCO-CE}, we describe how to use cross-validation, a
particular technique based on the bias-variance tradeoff, to modify
the CE method. Sec.~\ref{sec:results} then presents performance
comparisons between this modified version of the CE method and the
conventional CE method. These comparisons are on continuous,
multimodal, unconstrained optimziation problems. We show that on these
problems, using the modified version of the CE method can
significantly improve optimization performance of the CE method, and
never worsens performance.

\section{The Bias-Variance Tradeoff}
\label{sec:BPVDerivation}
Consider a given random variable $\targetrv$ and a random variable
that we can modify, $\estimatorrv$. We wish to use a sample $\estimatorvalue$ of $\estimatorrv$
as an estimate of a sample $\phi$ of $\targetrv$. The mean squared error
between such a pair of samples is a sum of four terms. The first term
reflects the statistical coupling between $\targetrv$ and $\estimatorrv$ and is conventionally ignored in bias-variance analysis. 
The second term reflects the inherent randomness in $\targetrv$ and is independent of the
estimator $\estimatorrv$. Accordingly, we cannot affect this term. In
contrast, the third and fourth terms depend on $\estimatorrv$. The third
term, called the bias, is independent of the precise samples of both
$\targetrv$ and $\estimatorrv$, and reflects the difference between the means
of $\targetrv$ and $\estimatorrv$. The fourth term, called the variance, is
independent of the precise sample of $\targetrv$, and reflects the
inherent randomness in the estimator as one samples it.  These last two
terms can be modified by changing the choice of the estimator. In
particular, on small sample sets, we can often decrease mean
squared error by introducing a small bias that causes a
large reduction the variance. While most commonly used in machine
learning, bias-variance tradeoffs are applicable in a much broader context and in a variety of
situations. 

\subsection{A Simple Derivation of the Bias-Variance Decomposition}
Suppose we have a Euclidean random variable $\targetrv$ taking
on values $\phi$ distributed according to a density function $p(\phi)$.  We
want to estimate a certain value $\phi$, that we cannot access directly, and that was obtained by sampling $p(\phi)$. We can, however, access a different
Euclidean random variable $\estimatorrv$ taking on values $\estimatorvalue$
distributed according to $p(\estimatorvalue)$. We can also modify the distribution of $\estimatorrv$, and we want to exploit the coupling between $\estimatorrv$ and $\targetrv$ to improve our estimate. Assuming a quadratic loss function, the quality of our estimate is measured by its Mean Squared
Error (MSE):
\begin{eqnarray}
\textrm{MSE}(\estimatorrv) &\equiv& \int p(\estimatorvalue, \phi) \,(\estimatorvalue - \phi)^2 \,d\estimatorvalue\, d\phi.
\label{eq:bpv}
\end{eqnarray} 
In standard bias-variance analysis, the statistical coupling of $\targetrv$ and $\estimatorrv$ is simply ignored without any justification\footnote{One could account for this coupling by using an additive correction term~\citep[see][]{wolp97}.},
and the distribution $p(\phi, \estimatorvalue)$ is replaced with the product of marginals, $p(\phi)p(\estimatorvalue)$. So our equation for MSE reduces to
\begin{eqnarray}
\textrm{MSE}(\estimatorrv) &=& \int p(\estimatorvalue) p(\phi) \,(\estimatorvalue - \phi)^2 \,d\estimatorvalue\,d\phi.
\label{eq:exp_cost}
\end{eqnarray} 

Using simple algebra, the right hand side of Eq.~\ref{eq:exp_cost} can
be written as the sum of three terms. The first is the variance of
$\targetrv$. This term is justifiably ignored in bias-variance analysis since it is beyond our control in designing the estimator $\estimatorrv$. The second term involves a mean that describes the
deterministic component of the error. This term depends on both the
distribution of $\targetrv$ and that of $\estimatorrv$, and quantifies how close the
means of those distributions are.  The third term is a
variance that describes stochastic variations from one sample to the
next. This term is independent of the random variable being
estimated. Formally, up to an overall additive constant, we can write
\begin{eqnarray}
\textrm{MSE}(\estimatorrv) &=& \int p(\estimatorvalue) (\estimatorvalue^2 - 2 \phi \estimatorvalue + \phi^2) d\estimatorvalue,\nonumber\\
&=& \underbrace{\int p(\estimatorvalue) \estimatorvalue^2 d\estimatorvalue} - 2 \phi \int p(\estimatorvalue) \estimatorvalue\, d\estimatorvalue + \phi^2,\nonumber\\ 
&=& \var(\estimatorvalue) + [\mean(\estimatorvalue)]^2 - 2 \phi \;\mean(\estimatorvalue) + \phi^2,\nonumber\\
&=& \var({\estimatorvalue}) + \underbrace{[\phi - \mean(\estimatorvalue)]^2}, \nonumber\\
&=& \textrm{variance} + \textrm{bias}^2.
\label{eq:biasPlusVariance}
\end{eqnarray}

In light of Eq.~\ref{eq:biasPlusVariance}, one way to try to reduce
MSE is to modify an estimator to trade bias
for variance. Some of the most well-known applications of such
bias-variance tradeoffs occur in parametric machine learning, where
many techniques have been developed to exploit that tradeoff. 
There are still some extensions of that tradeoff that could be applied in parametric machine learning that have been ignored by that
community. 
$ $


\section{Applications of the Bias-Variance Tradeoff}
\label{sec:applications}

In this section, we describe some applications of the bias-variance
tradeoff. First, we provide a simple, concrete example that elaborates how conventional bias-variance techniques ignore the statistical coupling described in the previous section. Then, we describe the application of bias-variance tradeoffs to Monte Carlo (MC) techniques for the
estimation of integrals, where, as for all unbiased estimators, the bias-variance tradeoff reduces to simple variance reduction. Next, we introduce the field of Monte
Carlo Optimization (MCO), and illustrate that there are subtleties involved that are absent in simple MC. Then, we describe the field of Parametric Machine Learning, which is mathematically
identical to MCO, and describe how they ignore some of the subtleties pertaining to MCO, but apply bias-variance analysis nonetheless.

\subsection{Supervised Learning}
Consider the simplest type of supervised machine learning problem, where the aim is to accurately predict the behavior of an input-output system, based on several `training examples'. In this example, we consider a finite input space $X$, and an output space $Y$ which is just the space of real numbers, and a deterministic input-output system, or
`target function' $f$ that maps each element of $X$ to a single
element of $Y$. To be precise, there is a `prior' probability density function
$p(f)$ over target functions, and it gets sampled to produce some
particular target function, $f$.  Next, $f$ is IID sampled at a set of
$m$ inputs to produce a `training set' ${\cal{D}}$ of input-output
pairs.

For simplicity of analysis, say we have a single fixed `prediction point' $x
\in X$.  Our  goal in supervised learning is to estimate
$f(x)$, but $f$ is not known. Accordingly, to perform the
estimation the training set is presented to a `learning algorithm',
which in response to the training set produces a guess $\hat{f}(x)$ for the
value $f(x)$.

This entire stochastic procedure defines a joint distribution $\pi(f,
{\cal{D}}, f(x), \hat{f}(x))$. In order to analyze the performance of the learning algorithm, we marginalize to obtain a distribution
$\pi(f(x), \hat{f}(x))$. Since $\hat{f}(x)$ is supposed to be an estimate of
$f(x)$, we can identify $\hat{f}(x)$ as the value $\estimatorvalue$ of the random
variable $\estimatorrv$ and $f(x)$ as the value $\phi$ of $\targetrv$. In other
words, we can define $p(\phi, \estimatorvalue) = \pi(f(x), \hat{f}(x))$. If we then compute the mean squared error in the estimate made by our learning algorithm for the value $f(x)$, we get Eq.~\ref{eq:bpv}. 

Of course, in general, $\targetrv$ and $\estimatorrv$, i.e., $f(x)$ and $\hat{f}(x)$, are statistically dependent; if they weren't, the learning algorithm gains nothing from knowing $\cal{D}$.
This dependence can be established by writing
\begin{eqnarray}
p(f(x), \hat{f}(x)) &=& \int d{\cal{D}} \; p(f(x), \hat{f}(x) \mid {\cal{D}}) \;
p({\cal{D}}), \nonumber\\
&=& \int d{\cal{D}} \; p(\hat{f}(x) \mid f(x), {\cal{D}}) \;p(f(x) \mid
{\cal{D}}) \; p({\cal{D}}),\nonumber \\
&=& \int d{\cal{D}} \; p(\hat{f}(x) \mid {\cal{D}}) \; p(f(x) \mid
{\cal{D}}) \; p({\cal{D}}).
\label{eq:BPVExact}
\end{eqnarray}
The first two steps comprise a straightforward application of Bayes' rule. In the third step, the conditioning on $\hat{f(x)}$ is removed because the guess of the learning algorithm is determined solely by the training set $\mathcal{D}$. Nevertheless, note that Eq.~\ref{eq:BPVExact} exact, and in general is \emph{not} the same as the product
\begin{eqnarray}
p(f(x)) \; p(\hat{f}(x)) &=& \left[\int d{\cal{D}} \; p(f(x) \mid {\cal{D}})
p({\cal{D}})\right]\left[\int d{\cal{D}} \; p(\hat{f}(x) \mid {\cal{D}})
p({\cal{D}})\right].
\label{eq:BPVApprox}
\end{eqnarray}
As we described in Sec.~\ref{sec:BPVDerivation}, conventional bias-variance analysis approximates Eq.~\ref{eq:BPVExact} by Eq.~\ref{eq:BPVApprox}.

\subsection{Monte Carlo Integration}

Monte Carlo methods are often the method of choice for estimating
difficult high-dimensional integrals. Consider a function $f \colon
X \rightarrow \mathbb{R}$, which we want to integrate over
some region $\mathcal{X} \subseteq X$, yielding the value
$F$, as given by
\begin{equation*}
F = \int_{\mathcal{X}} dx\, f(x).
\end{equation*} 
We can view $F$ as a degenerate random variable $\targetrv$, with density function given
by a Dirac delta function centered on $F$.  Therefore, the variance of $\targetrv$
is 0, and Eq.~\ref{eq:biasPlusVariance} is exact.

A popular MC method to estimate this integral is importance sampling
\cite[][]{roca04}. This exploits the law of large numbers as follows:
i.i.d. samples $x^{(i)}, \; i=1,\ldots, m$ are generated from a
so-called importance distribution $h(x)$ that we control, and the
associated values of the integrand, $f(x^{(i)})$ are computed. Denote
these `data' by 
\begin{equation}
\mathcal{D} = \{(x^{(i)}, f(x^{(i)}), \;i=1,\ldots,m\}.
\label{eq:data}
\end{equation} 
Now,
\begin{eqnarray*}
\phi\;\; = \;\;F &=& \int_{\mathcal{X}} dx\, h(x) \frac{f(x)}{h(x)},\nonumber\\
&=& \lim_{m \rightarrow \infty} \frac{1}{m}\sum_{i=1}^m \frac{f(x^{(i)})}{h(x^{(i)})}\textrm{ with probability 1.}
\end{eqnarray*}
Denote by $\estimatorrv$ the random variable with value given by the \emph{finite}
sample average for $\mathcal{D}$:
\begin{equation*}
\estimatorvalue = \frac{1}{m}\sum_{i=1}^m \frac{f(x^{(i)})}{h(x^{(i)})}.
\end{equation*}
We use $\estimatorrv$ as our statistical estimator for $\targetrv$, as we
broadly described in Sec.~\ref{sec:intro}.  Assuming a quadratic
loss function, $L(\estimatorvalue, \phi) = (\phi - \estimatorvalue)^2$, the bias-variance decomposition described
in Eq.~\ref{eq:biasPlusVariance} applies \emph{exactly}. 
It can be shown that the estimator $\estimatorrv$ is unbiased, that is,
$\mean{}\estimatorrv = \phi$, where the mean is over samples of $h$. Consequently, the MSE of this estimator is just its variance. The
choice of sampling distribution $h$ that minimizes this variance is
given by~\citep*[see][]{roca04}
\begin{equation*}
h^{\star}(x) = \frac{|f(x)|}{\int_{\mathcal{X}} |f(x')| dx'}.
\end{equation*}

By itself, this result is not very helpful, since the equation for the
optimal importance distribution contains a similar integral to the one
we are trying to estimate. For non-negative integrands $f(x)$, the
VEGAS algorithm \citep{lepa78} describes an adaptive method to find
successively better importance distributions, by iteratively
estimating $\targetrv$, and then using that estimate to generate the next
importance distribution $h$. In the case of this and other unbiased estimators,
there is no tradeoff between bias and variance, and minimizing MSE is
achieved by minimizing variance.

\subsection{Monte Carlo Optimization}
Instead of a \emph{fixed} integral to evaluate, consider a parametrized integral
\begin{equation*}
F(\theta) = \int_{\mathcal{X}} dx\,f_{\theta}(x).
\end{equation*}
Further, suppose we are interested in finding the value of the
parameter $\theta \in \Theta$ that minimizes $F(\theta)$: 
\begin{equation*}
\theta^{\star} = \arg\min_{\theta \in \Theta} F(\theta).
\end{equation*}
In the case where the functional form of $f_{\theta}$ is not explicitly known,
one approach to solve this problem is a technique called Monte Carlo
Optimization (MCO) \citep*[see][]{erno98}, involving repeated MC estimation of the integral in question with adaptive modification\footnote{The similarity to the CE method is quite clear.} of the parameter $\theta$. 

We proceed analogously to the preceding section. Whereas in MC, there was no parameter $\theta$ and we had a single Dirac-delta distribution, we now have a set of such distributions, corresponding to the values $\theta$. Accordingly, we introduce the $\theta$-indexed vector random variable $\targetrv$, each of whose components $\targetrv_{\theta}$ has a degenerate Dirac-delta distribution about the associated
value $F(\theta)$. 
Next, we introduce our estimator random-variable, a similar $\theta$-indexed vector random variable $\estimatorrv$ each of whose components $\estimatorrv_{\theta}$ can be sampled to estimate $\targetrv_{\theta}$. Regardless of how $\estimatorrv$ is defined, given a sample of $\estimatorvalue$, one way
to estimate $\theta^{\star}$ is
\begin{equation*}
\hat{ {{\theta}}} ^{\star} = \arg\min_{\theta \in \Theta} \estimatorvalue_{\theta}
\end{equation*}
We call this approach `natural' MCO. 

For example, let $\mathcal{D}$ be a data set as described in Eq.~\ref{eq:data}. Then for every $\theta$, any
sample of $\mathcal{D}$ provides a sample of $\estimatorrv_{\theta}$, which is the associated estimate 
\begin{equation*}
\estimatorvalue_{\theta}\;\; = \;\;\hat{F}(\theta) = \frac{1}{m}\sum_{i=1}^m \frac{f_\theta(x^{(i)})}{h(x^{(i)})},
\end{equation*}  
Given this choice for $\estimatorrv$, the `natural MCO' approach is to estimate the optimal $\theta$ as follows.
\begin{equation}
\hat{\theta}^{\star} =\arg\min_{\theta \in \Theta}\frac{1}{m}\sum_{i=1}^m \frac{f_{\theta}(x^{(i)})}{h(x^{(i)})}.
\label{eq:naiveMCO}
\end{equation}
As we shall see, this does not work well in practice, and we therefore call this importance-sampling application of natural MCO  `naive' MCO.

In general, consider \emph{any} algorithm that estimates $\theta^{\star}$ as a
single-valued function of ${\estimatorvalue}$. The estimate of $\theta^\star$ produced by that
algorithm is itself a random variable, since it is a function of the
random variable $\estimatorrv$. Call this random variable $\hat{\Theta}^{\star}$, taking on
values $\hat{ {\theta}}^{\star}$. Any MCO algorithm is defined by
$\hat{\Theta}^{\star}$; that random variable encapsulates the estimate made by the algorithm.

To analyze the error of such an algorithm, consider the associated
random variable $F(\hat{\Theta}^{\star})$, taking on the exact values of the integral $F(\hat{\theta}^{\star})$. Since our aim in MCO is to minimize $F(\theta)$, it is reasonable to propose that the cost of estimating $\theta^{\star}$ by $\hat{\theta}^{\star}$ is nothing but the difference between $F(\hat{\theta}^{\star})$ and the true minimal value of the integral, $F({\theta^{\star}}) = \min_\theta F(\theta)$. In other words, we adopt the loss function 
\begin{equation}
L(\hat{\theta}^{\star}, {\theta^{\star}})
\triangleq F(\hat{\theta}^{\star}) - F({\theta^{\star}}).
\end{equation}
This is in contrast to our discussion on MC integration, which involved quadratic loss. Up to an unknown additive constant, this loss function is nothing but 
$F(\hat{\theta}^{\star})$. This additive constant
$F({\theta^{\star}})$ is fixed by the MCO problem at
hand, is beyond our control, and is not affected by our algorithm. Consequently, we ignore that additive constant, and write out the associated expected loss:
\begin{eqnarray}
\mean(L) &=& \int d\hat{\theta}^{\star} p(\hat{\theta}^{\star}) F(\hat{\theta}^{\star}).
\label{eq:expectedLoss}
\end{eqnarray}
Now change coordinates in this integral from the values of the scalar
random variable $\hat{\theta}^{\star}$ to the values of the underlying
vector random variable ${\estimatorrv}$.  The expected loss now becomes
\begin{eqnarray*}
\mean(L) &=& \int d{\estimatorvalue} \;p({\estimatorvalue})
F({\hat{\theta}}^\star(\estimatorvalue)).
\end{eqnarray*}
The natural MCO algorithm provides some insight into these
results.  For that algorithm,
\begin{eqnarray}
\mean(L) &=& \int d{\estimatorvalue} \;p({\estimatorvalue})
F(\arg\min_{\theta}{\estimatorvalue}_{\theta}) \nonumber \\
&=& \int d{\estimatorvalue}_{\theta_1} d\estimatorvalue_{\theta_2} \ldots \;p({\estimatorvalue_{\theta_1}, \estimatorvalue_{\theta_2}, \ldots}) F(\arg\min_{\theta}{\estimatorvalue}_{\theta}).
\label{eq:natural_loss}
\end{eqnarray}
For any fixed $\theta$, there is an error between samples
$\estimatorvalue_{\theta}$ and the true value $F(\theta)$. Bias-variance
considerations apply to this error, exacty as in the discussion of MC
above. In MCO, however, we are not concerned with $\estimatorvalue$ for a single
component $\theta$, but rather for a set of $\theta$'s. 

The simplest such case is where the components $\estimatorvalue_{\theta}$ are independent. Even so, $\arg\min_\theta \estimatorvalue_{\theta}$ is distributed according
to the laws for extrema of multiple independent random
variables, and this distribution depends on higher-order moments of each random variable
$\estimatorvalue_{\theta}$. This means that $\mean(L)$ also depends on such higher-order
moments. Only the first two moments, however, arise in the bias and
variance for any single $\theta$. Thus, even in the simplest possible case, the bias-variance considerations for the individual $\theta$ do not provide a complete analysis.

In most cases, the components of $\estimatorvalue$ are $not$ independent. Therefore, in order to analyze $\mean(L)$, in addition to higher moments of the distribution for each individual $\theta$, we must now also consider higher-order moments coupling the estimates
${\estimatorvalue}_{\theta}$ for different $\theta$. Conventional bias-variance analysis for MCO is therefore incomplete on three fronts: ignoring the coupling between $\targetrv$ and $\estimatorrv$, ignoring higher order moments for individual $\theta$'s, and ignoring moments coupling different $\theta$'s.

Due to the coupling between $\theta$'s, it may be quite acceptable for  the individual
components ${\estimatorvalue}_{\theta}$ to have both a large bias and a large
variance, as long as the covariances are large. Large
covariances would ensure that if some ${\estimatorvalue}_{\theta}$ were incorrectly large, then $\estimatorvalue_{\theta'}$, for all $\theta' \neq \theta$ would also be incorrectly large. This would preserve the ordering of $\theta$'s under $\targetrv$. So, even with large bias and variance for each $\theta$, the estimator as a whole would still work well. If we could exploit this insight, it may be possible to come up with weaker requirements for accurate estimators. Nevertheless, we can ignore these insights, and impose a stronger requirement: design estimators
${\estimatorvalue}_{\theta}$ with sufficiently small bias plus variance for
each single $\theta$. More precisely, suppose that those terms are
very small on the scale of differences $F(\theta) - F(\theta')$ for
any $\theta$ and $\theta'$. Then, by Chebychev's inequality, we know that
the density functions of the random variables $\estimatorrv_{\theta}$ and
$\estimatorrv_{\theta'}$ have almost no overlap.  Accordingly, the probability that a sample $\estimatorvalue_{\theta} - \estimatorvalue_{\theta'}$ has the opposite sign of $F(\theta) -
F(\theta')$ is almost zero. 

Evidently, $\mean(L)$ is generally determined by a complicated
relationship involving bias, variance, covariance, and higher
moments. Natural MCO in general, and naive MCO in particular, ignore all of these effects, and consequently, 
often perform quite poorly in practice. In the next section we discuss some ways of
addressing this problem.

\subsection{Parametric Machine Learning}

There are many versions of the basic MCO problem described in the
previous section. Some of the best-explored arise in parametric
density estimation and parametric supervised learning, which together
comprise the field of Parametric machine Learning (PL). In particular, parametric supervised learning attempts to solve
\begin{equation*}
\arg\min_{\theta \in \Theta} \int dx\, p(x) \int dy \,p(y\mid
x)f_{\theta}(x).
\end{equation*}
Here, the values $x$ represent inputs, and the values $y$ represent
corresponding outputs, generated according to some stochastic process
defined by a set of conditional distributions $\{p(y\mid x),\;
x\in\mathcal{X}\}$. Typically, one tries to solve this problem
by casting it as a single-stage MCO problem, For instance, say we adopt a
quadratic loss between a predictor $z_{\theta}(x)$ and the true value
of $y$. Using MCO notation, we can express the associated
supervised learning problem as finding $\arg\min_\theta F(\theta)$,
where
\begin{eqnarray}
l_{\theta}(x) &=& \int dy\, p(y\mid x)\, (z_{\theta}(x) - y)
^2,\nonumber\\
f_{\theta}(x) &=& p(x)\, l_{\theta}(x),\nonumber\\
F(\theta) &=& \int dx\, f_{\theta}(x).
\label{eq:sup_learn_1}
\end{eqnarray}

Next, the argmin is estimated by minimizing a sample-based estimate of
the $F(\theta)$'s. More precisely, we are given a `training set' of
samples of $p(y \mid x)\,p(x)$, \{$(x^{(i)}, y^{i}) i = 1, \ldots,
m$\}. This training set provides a set of associated estimates of
$F(\theta)$:
\begin{equation*}
\hat{F}(\theta) = \frac{1}{m}\sum_{i=1}^m l_{\theta}(x^{(i)}).
\end{equation*} 
These are used to estimate $\arg\min_\theta F(\theta)$, exactly as in
MCO.  In particular, one could estimate the minimizer of $F(\theta)$
by finding the minimium of $\hat{F}(\theta)$, just as in natural
MCO. As mentioned above, this MCO algorithm can perform very poorly in
practice. In PL, this poor performance is called `overfitting the
data'.

There are several formal approaches that have been explored in PL to
try to address this `overfitting the data'. Interestingly, none
are based on direct consideration of the random variable
${\hat{\theta}}^\star(\estimatorrv)$ and the ramifications of its
distribution for expected loss (cf.  Eq.~\ref{eq:natural_loss}). In
particular, no work has applied the mathematics of extrema of multiple
random variables to analyze the bias-variance-covariance tradeoffs
encapsulated in Eq.~\ref{eq:natural_loss}.

The PL approach that perhaps comes closest to such direct
consideration of the distribution of
$F(\hat{\theta}^\star)$ is uniform convergence theory,
which is a central part of Computational Learning Theory
\citep[see][]{angl92}. Uniform convergence theory starts by crudely
encapsulating the quadratic loss formula for expected loss under
natural MCO, Eq.~\ref{eq:natural_loss}. It does this by considering
the worst-case bound, over possible $p(x)$ and $p(y \mid x)$, of the
probability that $F(\theta^\star)$ exceeds $\min_\theta
F(\theta)$ by more than $\kappa$. It then examines how that bound
varies with $\kappa$. In particular, it relates such variation to
characteristics of the set of functions $\{f_\theta : \theta \in
\Theta\}$, e.g., the `VC dimension' of that set~\citep[see][]{vapn82,vapn95}.

Another approach is to apply bias-plus-variance
considerations to the $entire$ PL algorithm
$\hat{\Theta}^\star$, rather than to each
$\estimatorrv_{\theta}$ separately. This approach is applicable for
algorithms that do not use natural MCO, and even for non-parametric
supervised learning. As formulated for parameteric supervised
learning, this approach combines the formulas in
Eq.~\ref{eq:sup_learn_1} to write
\begin{eqnarray*}
F(\theta) &=& \int dx\, dy \; p(x) p(y \mid x) (z_\theta(x) - y)^2.
\end{eqnarray*}
This is then substituted into Eq.~\ref{eq:expectedLoss}, giving
\begin{eqnarray}
{\mean(L)} &=& \int d{\hat{\theta}}^\star  dx \,dy \; p(x)\, p(y \mid x)\,
p({\hat{\theta}}^\star) (z_{{{\hat{\theta}}}^\star}(x) - y)^2
\nonumber \\
&=& \int dx \; p(x) \left[ \int  d{\hat{\theta}}^\star dy \; p(y
\mid x)
p({\hat{\theta}}^\star) (z_{{{\hat{\theta}}}^\star}(x) - y)^2 \right].
\label{eq:sup_learn_bpv}
\end{eqnarray}
The term in square brackets is an $x$-parameterized
expected quadratic loss, which can be
decomposed into a bias, variance, etc., in the usual way. This formulation eliminates any direct concern for issues like
the distribution of extrema of multiple random variables, covariances
between $\estimatorrv_{\theta}$ and $\estimatorrv_{\theta'}$ for different
values of $\theta$, and so on. In some ways, this imposes stronger conditions on the estimators we design, and consideration of the issues outlined above may give weaker conditions on how to increase estimation accuracy.

There are numerous other approaches for addressing the problems of
natural MCO that have been explored in PL. Particulary important among
these are Bayesian approaches, e.g.,
~\cite{buwe91,berg85,mack03}. Based on these approaches, as well as on intuition, many powerful techniques for addressing
data-overfitting have been explored in PL, including regularization,
cross-validation, stacking, bagging, etc. Essentially all
of these techniques can be applied to $any$ MCO problem, not just PL
problems. Since many of these techniques can be justified using
Eq.~\ref{eq:sup_learn_bpv}, they provide a way to exploit the
bias-variance tradeoff in other domains besides PL.

\section{PLMCO-CE}
\label{sec:PLMCO-CE}

In this section, we first present a review of the CE method for continuous problems, and describe various components in the terms used thus far. 

\citet{krpo06} describe the application of the CE method to continuous multi-extremal problems. They describe the use parametrized continuous distributions, in particular, multivariate Gaussian and multivariate Gaussian mixtures, to perform adaptive importance sampling for optimization of continuous functions. The problem at hand, in our notation, is as follows. Find
\begin{equation}
\textrm{minimize}_{x \in \mathcal{X}}\; G(x).
\end{equation}
This is then converted to an Associated Stochastic Problem (ASP) over $\theta$-parametrized probability distributions $q_{\theta}$ over $\mathcal{X}$. Let $\gamma^{\star} = \min_{x\in\mathcal{X}}[G(x)]$, and let $I_{\{\cdot\}}$ be the indicator function. We want to
\begin{equation}
\textrm{maximize}_{\theta} \;\mean_{q_{\theta}} I_{\{G(x) \leq \gamma^{\star}\}}
\end{equation}
This solution to this problem is a degenerate Dirac-delta function centered on the optimal $x$, provided that the set of allowable $\theta$ permits parametrization of such degenerate distributions. Note also that sampling this degenerate distribution gives the answer to the original problem, so in a way, the two problems are equivalent. 

The ASP is actually solved using a homotopy method, where one solves a sequence of optimization problems with progressively decreasing value of $\gamma$. The algorithm proceeds as follows. $\theta$ is initialized to some $\theta_0$. At each iteration $t \geq 0$, a set of samples is drawn from $q_{\theta_t}$, and $\gamma_{t+1}$ is chosen to correspond to the the best $\kappa$ percentile of these samples, and and $\theta_{t+1}$ is chosen so as to minimize a KL divergence (or cross entropy) from the suitably normalized indicator distribution (also called a Heaviside distribution) given by
\begin{equation}
p_{\gamma}(x) \propto \Theta_{\gamma}(x) = \left\{
\begin{array}{ll}
1,&G(x)\leq \gamma,\\
0,&\textrm{otherwise.} 
\end{array}
\right.
\end{equation}

The problem of computing $\theta_{t+1}$ is actually an MCO problem: we use a set of samples to search for the $\theta$ that minimizes a parametrized integral
\begin{equation}
\theta^{\star}_{t+1} = \arg\min_{\theta} \mathrm{KL}(p_{\gamma} \| q_{\theta}) = -\arg\min_{\theta}\int_{x \in\mathcal{X}} dx\, p_{\gamma}(x) \ln\left(q_{\theta}(x)\right).
\end{equation}
If we keep track of the $q_{\theta}$ that each of the elite samples was generated from, we can then properly use importance sampling to estimate the above integral. Also, the normalization constant turns out to be irrelevant, enabling the use of $\Theta_{\gamma}$ rather than $p_{\gamma}$. The importance-sampled estimate of of interest is therefore
\begin{equation}
\sum_{i=1}^m \dfrac{\Theta_{\gamma}(\xith{x}{i})}{q_{\theta_{k_i}}(\xith{x}{i})} \ln\left(q_{\theta}(\xith{x}{i})\right)
\end{equation}
Here, $q_{\theta_{k_i}}$ is the actual parametrized distribution that generated the $i^{\rm th}$ sample. In the CE method, however, the denominator of the likelihood ratio is for some reason ignored, and the importance ratio is just set to unity for the elite samples, and zero for all the others. 

\subsection{Cross-validation for the CE Method}
Using the notation from Sec.~\ref{sec:applications}, we see that the CE method actually performs naive MCO: it uses importance sampling to generate a finite-sample estimate of a parametrized integral, then estimates the integral-optimizing parameters by simply minimizing the finite-sample sum. It is known from extensive experience with learning algorithms in the PL community that such an approach is bound to perform poorly with small sample sizes. In other words, it will suffer from large errors caused by large variance, stemming from overfitting the data. Indeed, to prevent such overfitting, \citet{krpo06} recommend the use of `dynamic smoothing' to prevent `premature shrinking' of the distributions. Also, the percentile $\kappa$ must be set appropriately: if it is too large, convergence to an optimum will be very slow, and if it is too small, the algorithm will converge prematurely, either to a local optimum, or worse, to a point that is not even a local minimum, and make extremely slow progress towards a minimum. 

The percentile $\kappa$, or equivalently, the elite size $N^{\rm elite}$, can be thought of as a hyperparameter: a data-independent parameter that affects algorithm performance. In a PL context, hyperparameters are often set using the technique of cross-validation, which works as follows: the given data is partitioned into two parts, a training set and a `held-out' test set. The learning algorithm (also called a training algorithm) is given the training data, and the estimated optimal parameters $\hat{\theta}^{\star}$ it generates are tested by estimating the value of the parametrized integral using the test data. The hyperparameters are then chosen so as to optimize this `held-out' performance. This can be done many ways, helpfully named 70-30 cross-validation, $k$-fold cross-validation, or leave-one-out cross-validation. As expected, in 70-30 cross-validation, the given data is randomly partitioned into two parts containing 70\% and 30\% of the data. The larger partition is used for training, and the smaller one for testing. In $k$-fold cross-validation, the data is divided into $k$ equal-sized partitions. $k-1$ of these are used for training, and the last one is used for testing. This procedure is repeated $k$ times, so that there have been $k$ estimates generated and tested; the hyperparameters are chosen to optimize average held-out performance. In leave-one-out cross-validation, all but one data point are used for training, and the left-out point is used for testing. Of course, leave-one-out cross-validation only works if the performance measure makes sense on a single point. In some cases, such as prediction error, it does. In this paper, we use 4-fold cross-validation to trade bias and variance of the CE `training algorithm' to dynamically pick the percentile value $\kappa_t$ for elite samples. For our experiments, we used a fixed set of parameters to implement dynamic smoothing, though these too could be chosen using cross-validation. 

In the case of mixture distributions, \citet{krpo06} use a technique called augmentation, where each elite sample is presumed to have come from a particular component of the mixture. They don't specify how exactly they perform the partitioning of elite samples. Instead, we use the well-known Expectation Maximization (EM) algorithm that is widely used for precisely this purpose: maximizing likelihood (or cross-entropy) using mixture distributions. Broadly speaking, instead of deterministically partitioning the elite samples, the EM algorithm iteratively `soft-partitions' the elite samples using probabilities over a latent variable to specify which mixture component they arose from. The number of mixture components specifies what's called a model class. In addition to picking hyperparameters, cross-validation is also used to select model classes, and this process is called model selection. We use cross-validation to adaptively pick the number of mixing components in the Gaussian mixture $q_{\theta_t}$.

In our case, we want to use cross-validation to pick $\gamma$, or equivalently, $\kappa$. Therefore, rather than use $\mean_{q_{\theta}}I_{\{G(x) \leq \gamma\}}$ as the parametrized integral to optimize, we choose $\mean_{q_{\theta}}[G(x)]$. Even though we use a different ASP, note that the resulting optimization problem is formally equivalent to the CE ASP: they share the same optimum point $\theta$, but differ in optimum value. We do, however, still use the CE method as our training algorithm. In other words, we use varying $\kappa$ in the CE algorithm to generate estimates $\hat{\theta}^{\star}$ of the optimal $\theta$, and use cross-validation to choose that value $\kappa$ that minimizes the importance-sampled estimate of $\mean_{q_{\theta}}[G(x)]$ on the held-out data. In other words, pick $\kappa$ to minimize
\begin{equation}
\sum_{i=1}^{m_{ho}} \dfrac{q_{\theta}(\xith{x}{i})}{q_{\theta_{k_i}}(\xith{x}{i})} G(\xith{x}{i}).
\end{equation}
Note that the sum is over the held out samples, denoted by the subscript $ho$ on the summation limit. Also note that we do keep track of the actual likelihood for each sample generated, and use the correct likelihood ratio in our importance-sampled estimate for held-out performance. 

We call this CE algorithm, where hyperparameters and models are adaptively picked using cross-validation, PLMCO-CE, since PL techniques are applied to improve MCO performance of the CE algorithm.

\section{Results and Discussion}
\label{sec:results}
In this section, we present performance comparisons between the conventional CE method and the new version presented in Sec.~\ref{sec:PLMCO-CE}, on a few simple test problems for unconstrained optimization.

\subsection{Test Problems}
The test problems are all analytic, low-dimensional, multi-extremal functions. The problems varied in dimensionality from 4 to 8, and each of them has properties that make it difficult for local optimizers to find the global optima: some, such as the Woods and Rosenbrock problems, are badly-scaled, others have multiple local minima, with the worst minima having the largest basin of attraction for gradient-based algorithms, and so on. A few of these problems have been used for testing continuous multi-extremal optimization using the CE method \citep{krpo06}. The 4-dimensional problems are the $n$-dimensional (extended) Rosenbrock, the Woods function, and the Shekel family of functions. The Hougen function is 5-dimensional. We also experimented with the 6-dimensional Hartman function and the 8-dimensional Rosenbrock. The $n$-dimensional extended Rosenbrock is defined as
\begin{equation}
G(x) = \sum_{i=1}^{n-1} \left[(1-x_i)^2 + 100(x_i^2 - x_{i+1})^2\right]
\end{equation}
The Woods function is defined as
\begin{eqnarray}
G(x) &=& 100(x_2 - x_1)^2 + (1 - x_1)^2 + 90(x_4 - x_3^2)^2 + (1 - x_3)^2  \nonumber\\
     &&+ 10.1[(1 - x_2)^2 + (1 - x_4)^2] + 19.8(1-x_2)(1 - x_4),
\end{eqnarray}
The Hougen function is described by \citet{krpo06}. The others are described in detail by \citet{disz78}. 

\subsection{Algorithms Tested}
We versions of CE algorithm described by \citet{krpo06}, by varying value of the elite size $\nelite$ and the number of components in the Gaussian mixture. We tested three values of $\nelite$ as a fraction of the number of samples taken: 5, 10, and 15\%. For each of these, we experimented with a single Gaussian and a 3-component Gaussian mixture. In the attached plots, the single Gaussian algorithms are denoted by CESxx, where xx represents $\nelite$ as a percentage of the population size, and the mixture-based algorithms are denoted by CEMxx. We also tested two PLMCO-CE algorithms, which dynamically chose a good value of $\nelite$ at each iteration. These algorithms used cross-validation to pick the value of $\nelite$ that optimized estimated held-out performance, i.e., $\mean{q_{\theta}} [G(x)]$, as described in the preceding section. The version using Gaussian mixtures also dynamically chose the number of mixing components (between 1 and 3) using cross-validation. For a given number of mixing components, the algorithm used cross-validation to find the optimal value of $\nelite$. This follows the rather standard machine-learning practice of optimizing hyperparameters for each model type using cross-validation, and then picking the model type that optimizes average held-out performance. The PLMCO-CE algorithms using single Gaussian and mixtures are denoted by CESX and CEMX respectively.

\subsection{Comparing Algorithm Performance}
We compare the performance of these algorithms by running each of them 100 times on all 8 problems. For a given trial of a test problem, the same set of initial samples was used for all algorithms. Therefore, performance differences between algorithms reflect directly on the choices made by the algorithms themselves. Of course, this common initial set was varied from trial to trial.

We used two metrics to analyze the results of these experiments. The first is mean performance, but since we do not have the true mean, we plot the sample mean of these 100 runs, accompanied by the associated 95\% confidence intervals shown as shaded regions. This, however, is not a comprehensive summary of algorithm performance: we also present a semilog plot of algorithm performance, where we plot $(G_{\rm best} - G^{\star})$ on a log scale against the number of function calls. This is standard practice in conventional optimization, but is seldom done in evolutionary optimization community. One reason to use a semilog scale is to be able to visually distinguish smaller values when displayed on the same plot as larger ones: often, we want to know if some algorithm consistently finds numbers 0.001 smaller than some other algorithm, but the numbers themselves are much larger than 0.001. In conventional optimization, the other major reason to use a semilog scale is to assess convergence rate. On these semilog plots, we cannot properly display confidence intervals for two reasons: the upper and lower confidence intervals will not be of the same size owing the logarithmic nature of the plot scale, and besides, the value of the mean minus the confidence interval may be negative, and this cannot be shown on a log-scale. Owing to the positivity of the variable $(G_{\rm best} - G^{\star})$, this random variable has a rather skewed distribution, and in such cases, it is often more informative to consider the median, since the mean is dominated by the large values. So, our second metric is median performance, and this is shown plotted on a semilog scale, along with best and worst performance. Of course, these extremal events are those that occurred on our 100 runs, and we cannot make rigorous claims regarding the true best-case or worst-case performance. 

\subsection{Results}
\begin{figure}
\begin{center}
\begin{tabular}{cc}
\includegraphics[width=3in]{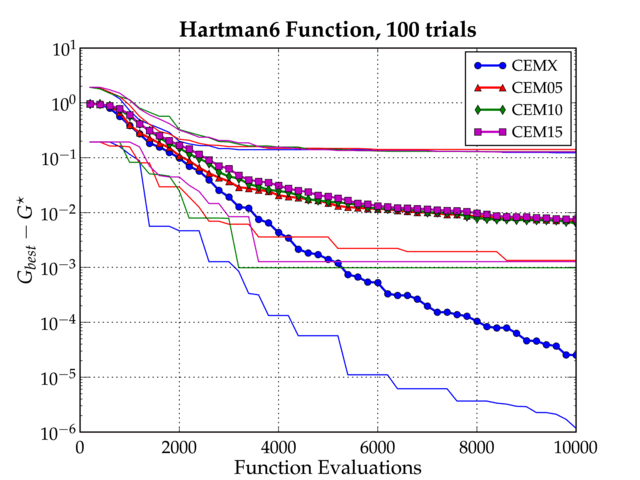}&
\includegraphics[width=3in]{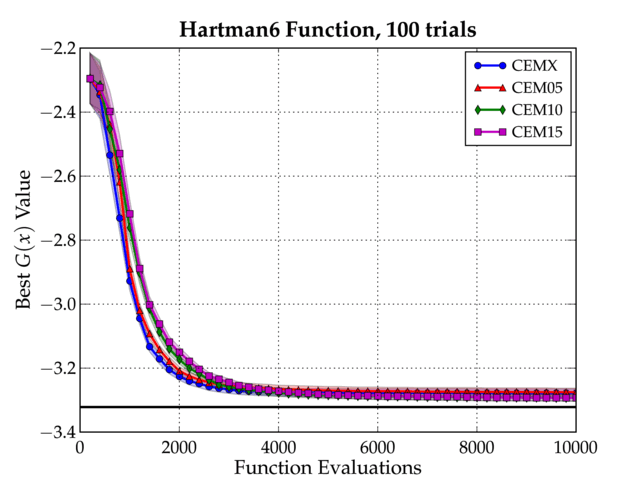}\\
a. Median performance, Gaussian mixtures. & b. Mean performance, Gaussian mixtures.\\\vspace{0.25in}\\
\includegraphics[width=3in]{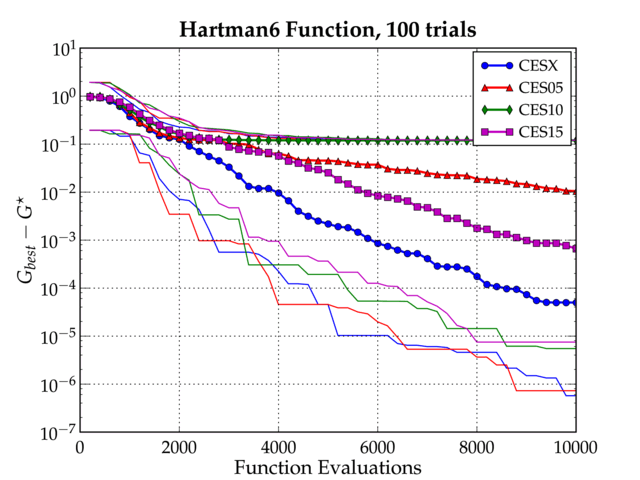}&
\includegraphics[width=3in]{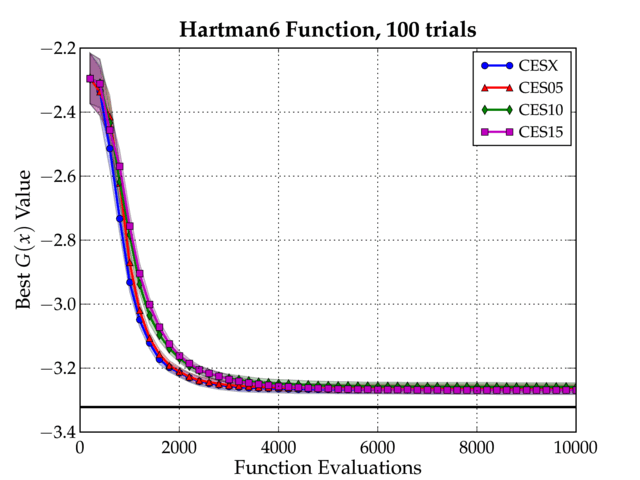}\\
c. Median performance, single Gaussian. & d. Mean performance, single Gaussian.
\end{tabular}
\caption{Performance comparison on Hartman6 function.}
\label{fig:performanceHartman6}
\end{center}
\end{figure}

We now present the performance comparisons discussed above. Fig.~\ref{fig:performanceHartman6} summarizes performance on the 6-dimensional Hartman problem. As we can see, the plot of mean performance does not show any discernible difference, but this is part of the reason for showing the semilog plot. We see that the median performance of the PLMCO-CE algorithm is many orders of magnitude better, both using Gaussian mixtures and a single Gaussian. This is not surprising on a multimodal problem, where `overfitting' to a set of unlucky data would lead to convergence to a local optimum. 

\begin{figure}
\begin{center}
\begin{tabular}{cc}
\includegraphics[width=3in]{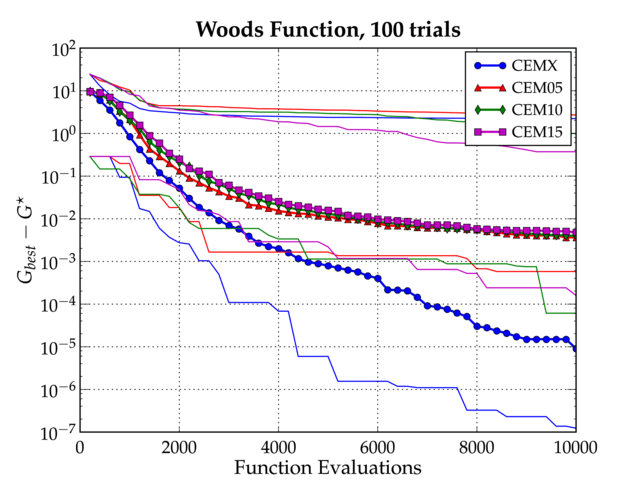}&
\includegraphics[width=3in]{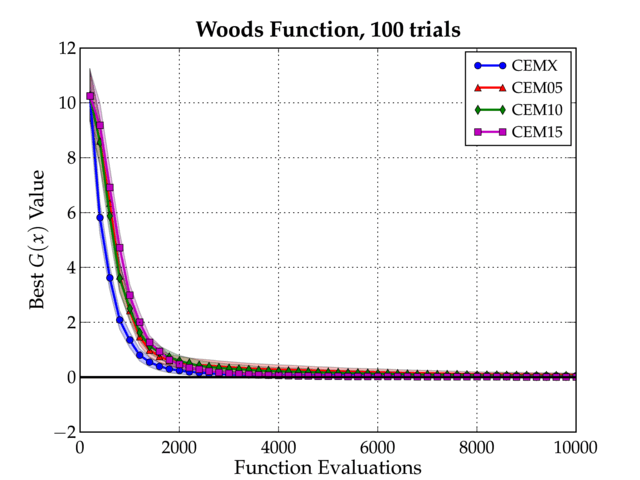}\\
a. Median performance, Gaussian mixtures. & b. Mean performance, Gaussian mixtures.\\\vspace{0.25in}\\
\includegraphics[width=3in]{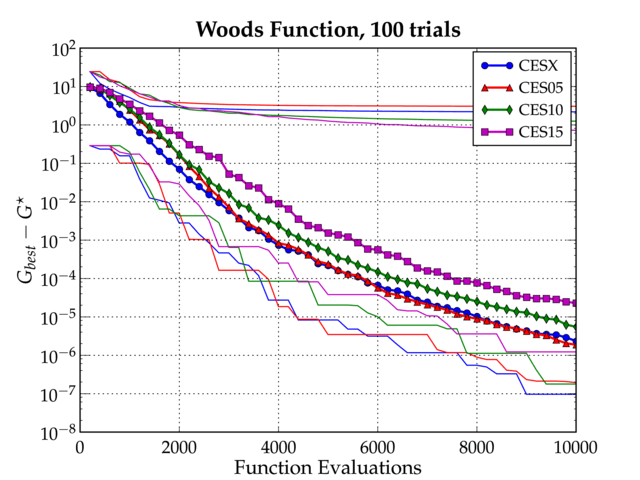}&
\includegraphics[width=3in]{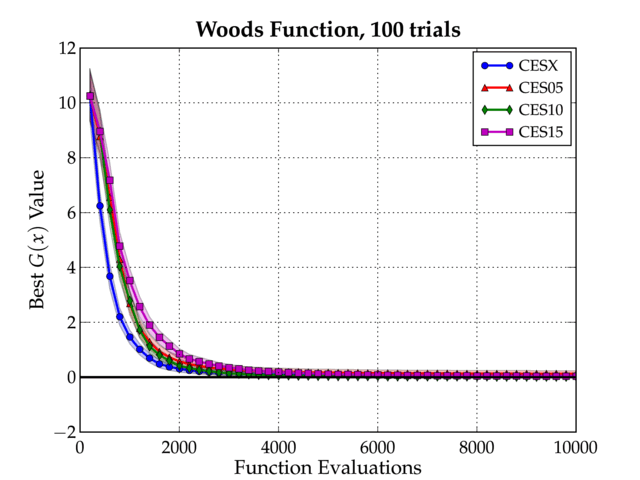}\\
c. Median performance, single Gaussian. & d. Mean performance, single Gaussian.
\end{tabular}
\caption{Performance comparison on Woods function.}
\label{fig:performanceWoods}
\end{center}
\end{figure}

\begin{figure}
\begin{center}
\begin{tabular}{ccc}
\includegraphics[width=2in]{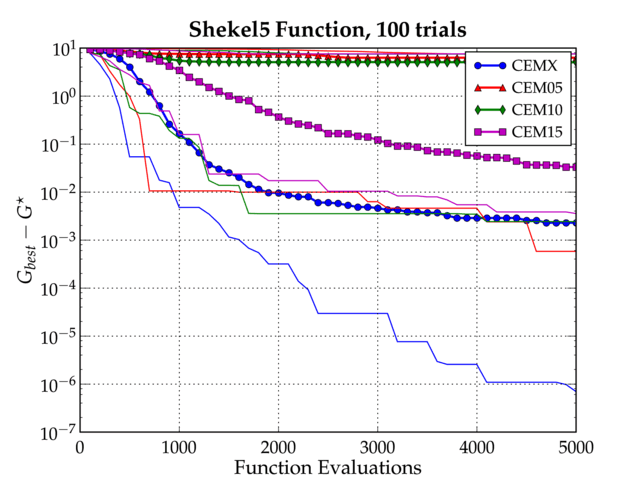}&
\includegraphics[width=2in]{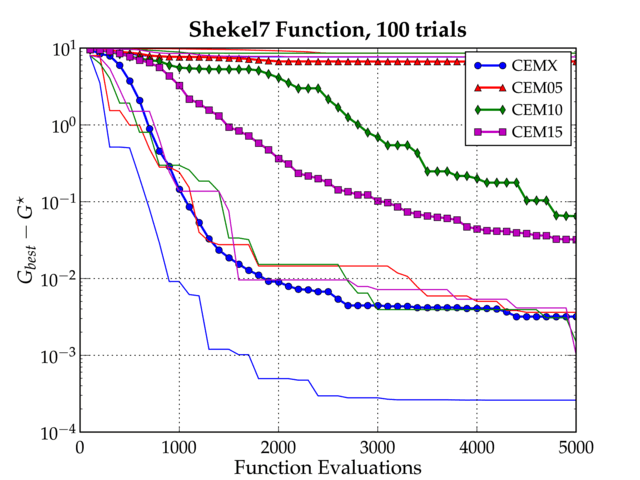}&
\includegraphics[width=2in]{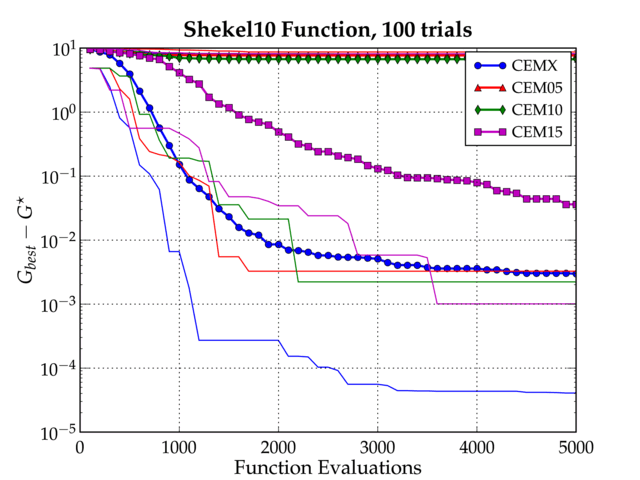}\\
\multicolumn{3}{c}{Median performance, Gaussian mixtures.}\\
\vspace{0.15in}\\
\includegraphics[width=2in]{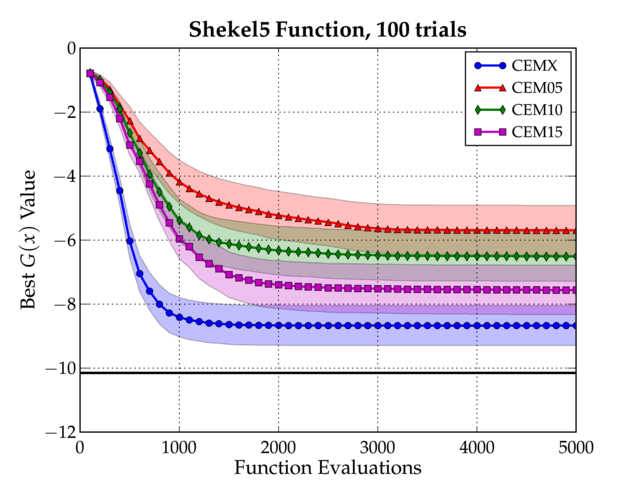}&
\includegraphics[width=2in]{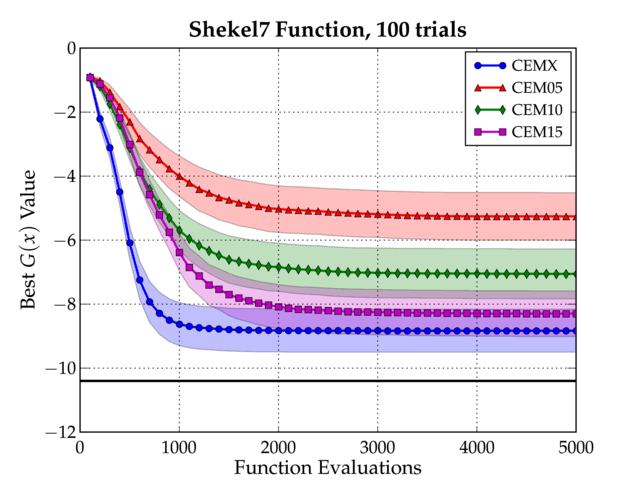}&
\includegraphics[width=2in]{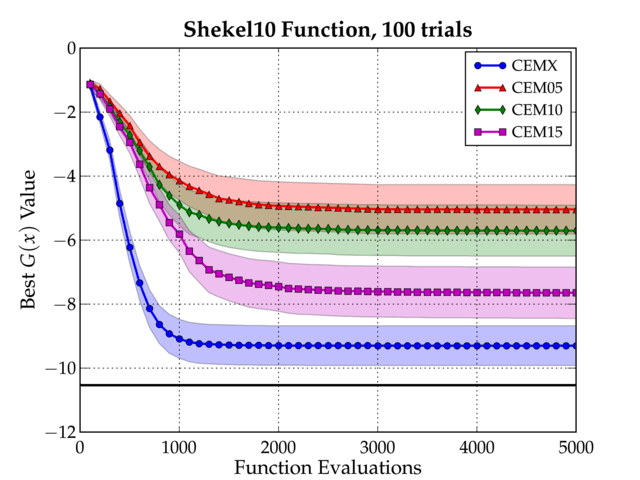}\\
\multicolumn{3}{c}{Mean performance, Gaussian mixtures.}\\
\vspace{0.15in}\\
\includegraphics[width=2in]{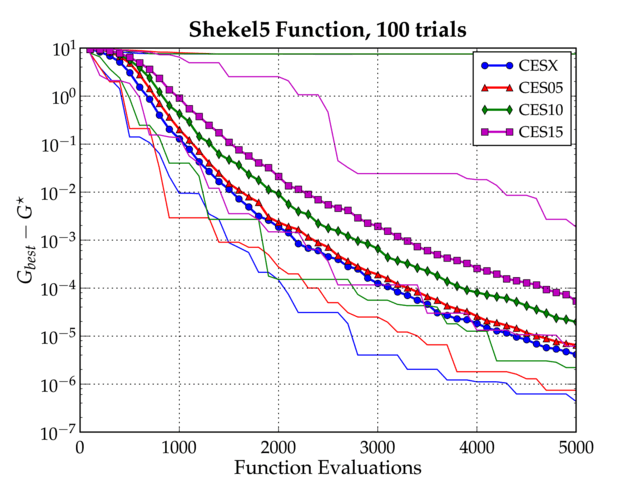}&
\includegraphics[width=2in]{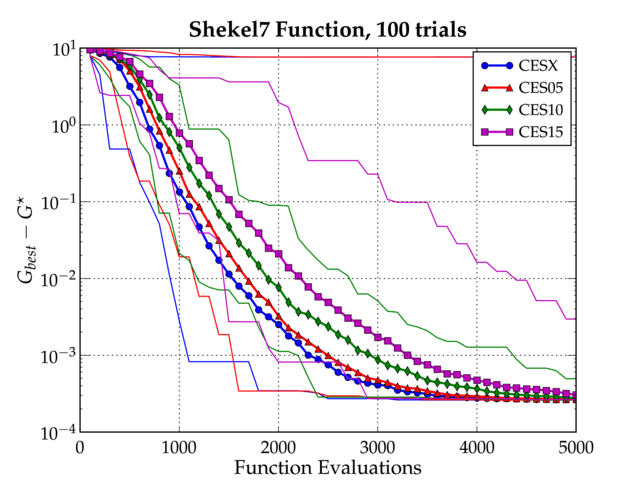}&
\includegraphics[width=2in]{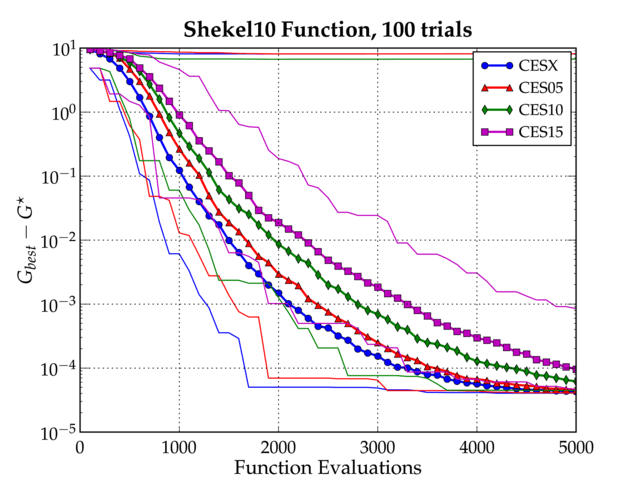}\\
\multicolumn{3}{c}{Median performance, single Gaussian.}\\
\vspace{0.15in}\\
\includegraphics[width=2in]{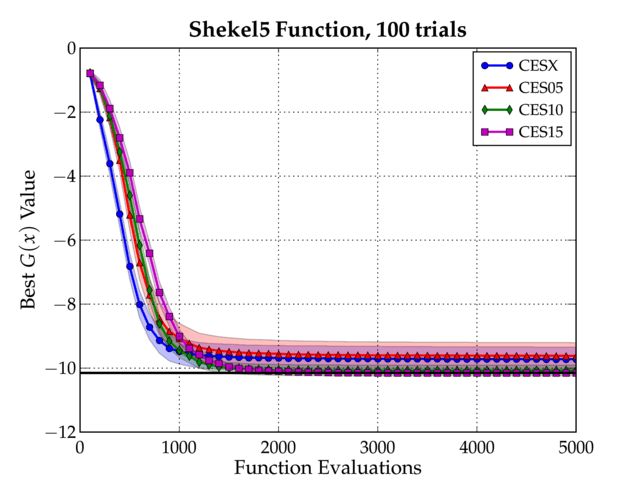}&
\includegraphics[width=2in]{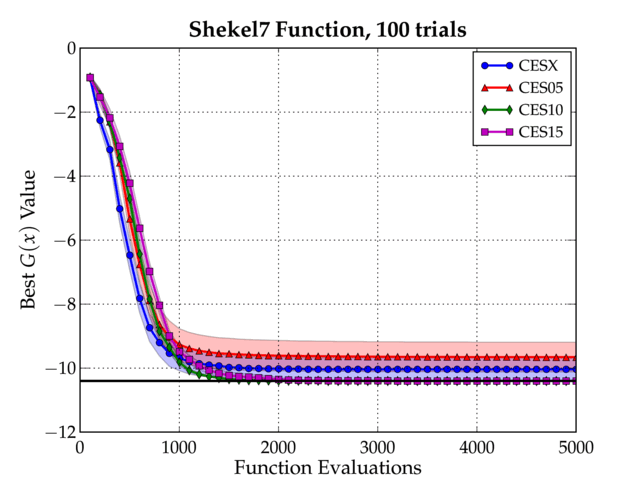}&
\includegraphics[width=2in]{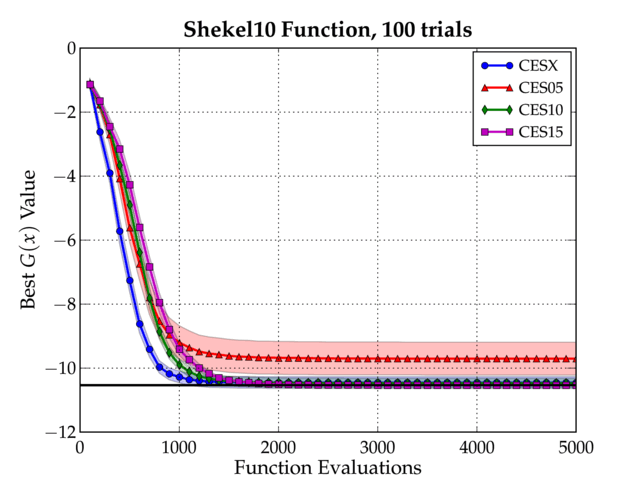}\\
\multicolumn{3}{c}{Mean performance, single Gaussian.}
\end{tabular}
\caption{Performance comparison on Shekel family of functions.}
\label{fig:performanceShekel}
\end{center}
\end{figure}

\begin{figure}
\begin{center}
\begin{tabular}{cc}
\includegraphics[width=3in]{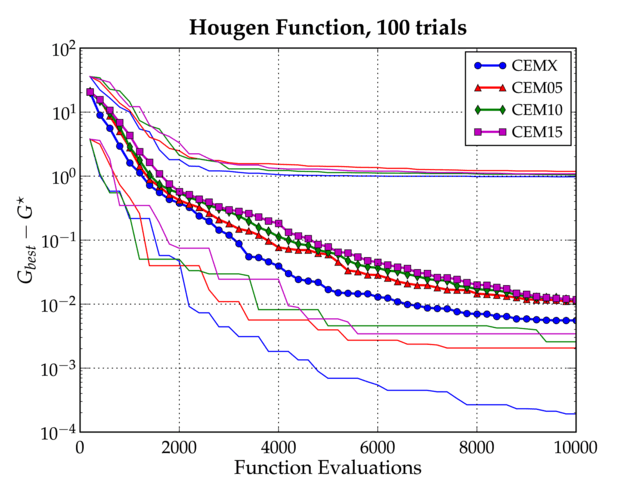}&
\includegraphics[width=3in]{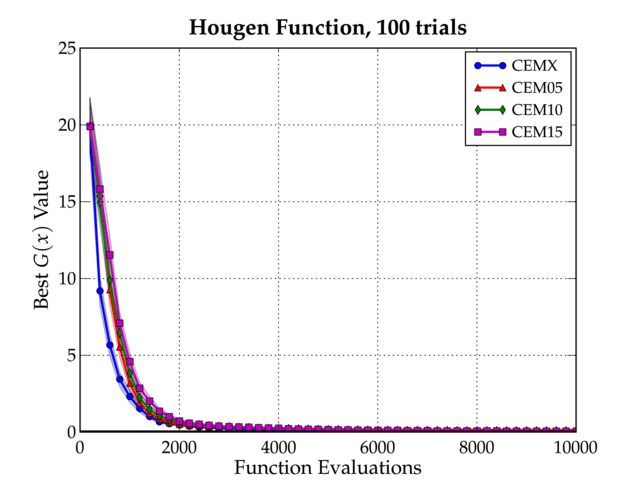}\\
a. Median performance, Gaussian mixtures. & b. Mean performance, Gaussian mixtures.\\\vspace{0.25in}\\
\includegraphics[width=3in]{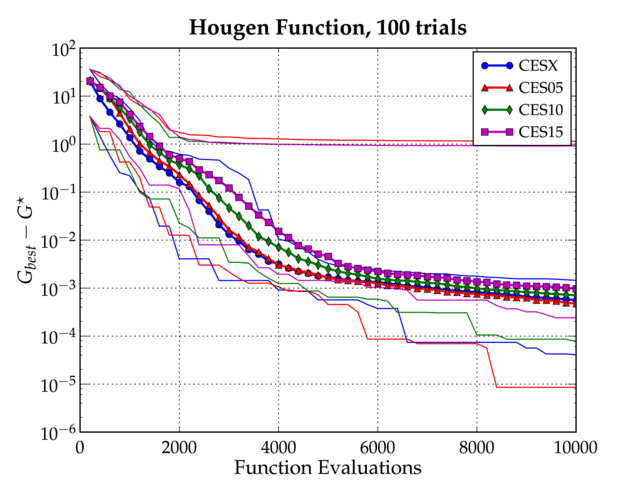}&
\includegraphics[width=3in]{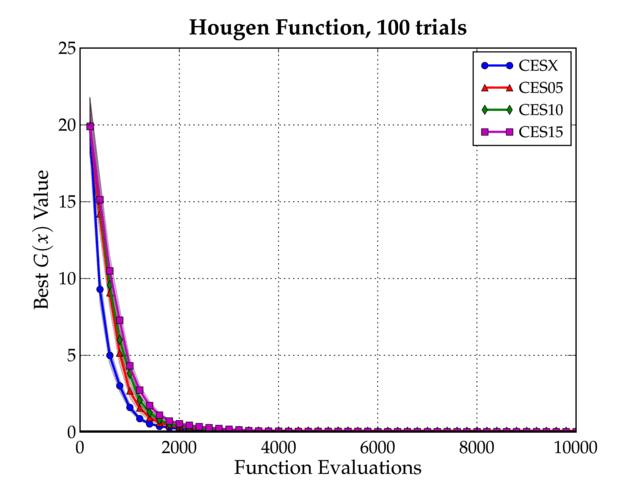}\\
c. Median performance, single Gaussian. & d. Mean performance, single Gaussian.
\end{tabular}
\caption{Performance comparison on Hougen function.}
\label{fig:performanceHougen}
\end{center}
\end{figure}

\begin{figure}
\begin{center}
\begin{tabular}{cc}
\includegraphics[width=2in]{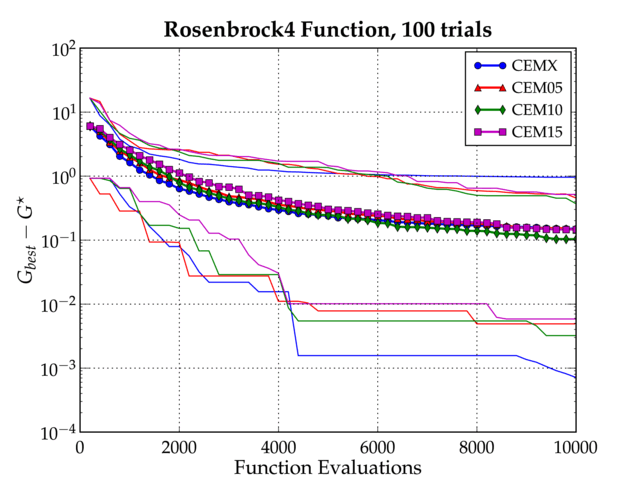}&
\includegraphics[width=2in]{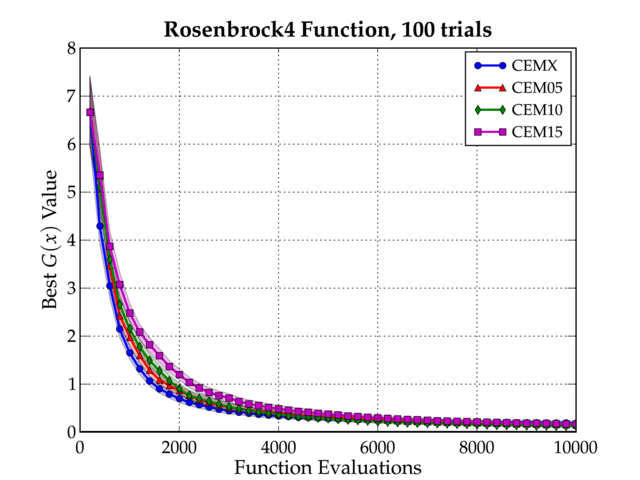}\\
a. Median performance, Gaussian mixtures. & b. Mean performance, Gaussian mixtures.\\\vspace{0.15in}\\
\includegraphics[width=2in]{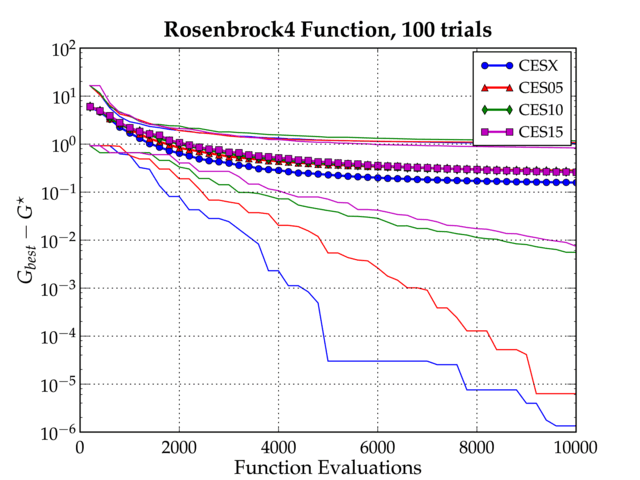}&
\includegraphics[width=2in]{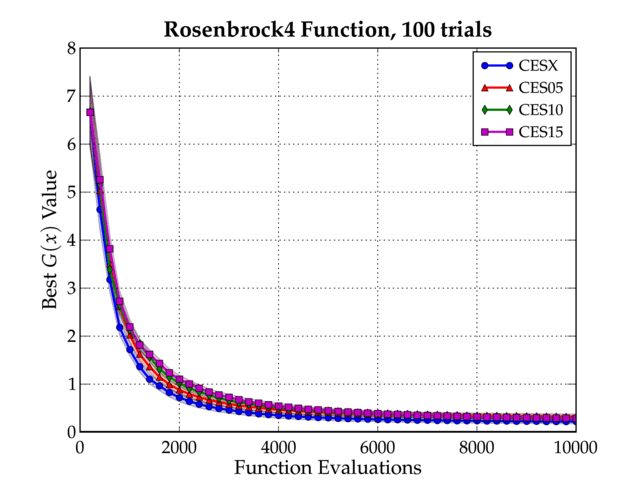}\\
c. Median performance, single Gaussian. & d. Mean performance, single Gaussian.\\
\vspace{0.15in}\\
\includegraphics[width=2in]{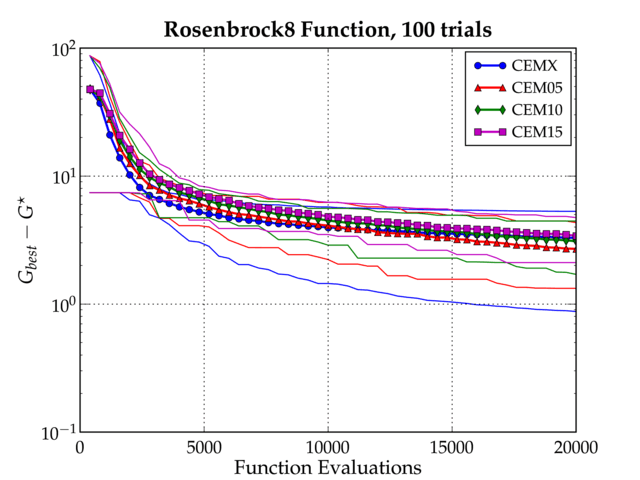}&
\includegraphics[width=2in]{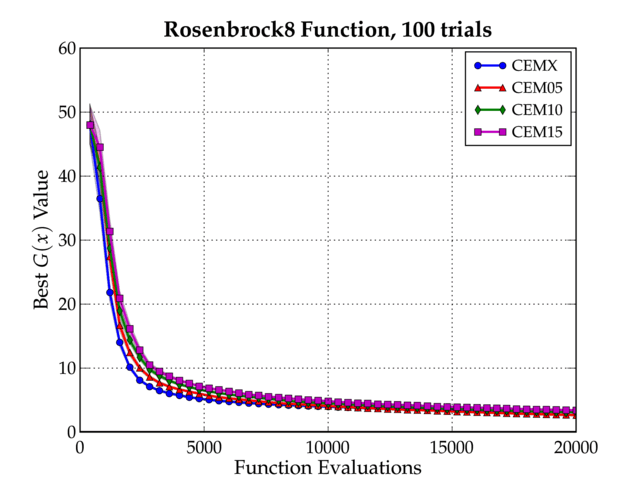}\\
e. Median performance, Gaussian mixtures. & f. Mean performance, Gaussian mixtures.\\\vspace{0.15in}\\
\includegraphics[width=2in]{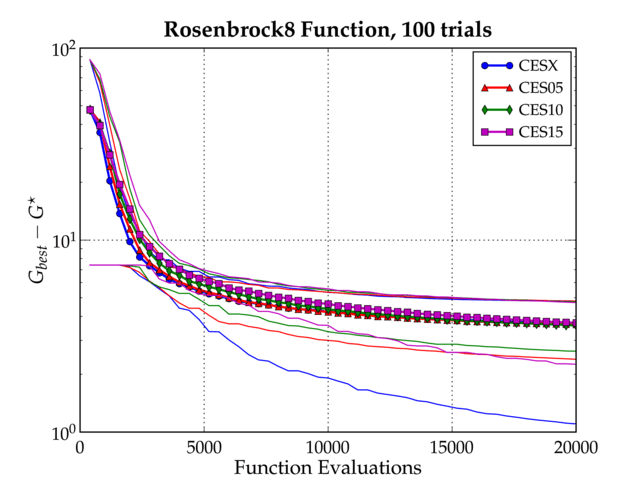}&
\includegraphics[width=2in]{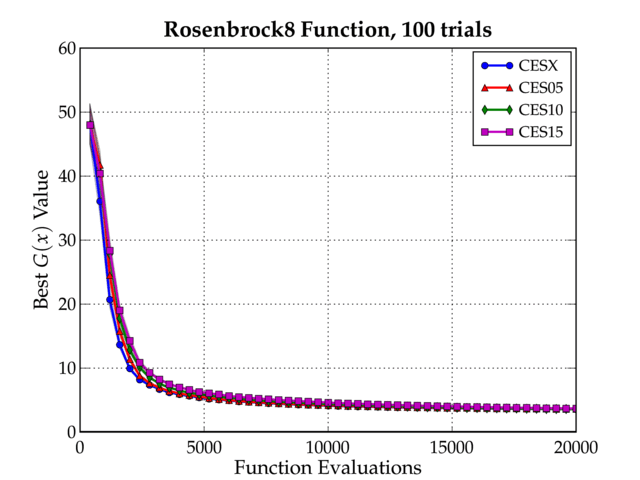}\\
g. Median performance, single Gaussian. & h. Mean performance, single Gaussian.
\end{tabular}
\caption{Performance comparison on $n$-dimensional Rosenbrock.}
\label{fig:performanceRosenbrock}
\end{center}
\end{figure}

Fig.~\ref{fig:performanceWoods} shows performance on the 4-dimensional Woods problem, which is unimodal, but is poorly scaled and has large flat regions. Using a single Gaussian, the mean performance of PLMCO-CE shows significant early gains, and yet final performance is no worse than the standard CE algorithm. This just means that PLMCO-CE quickly converges to the optimum. With mixtures, median performance, as usual, is greatly improved. 

Fig.~\ref{fig:performanceShekel} compares results on the Shekel family of functions, each of whose members is 4-dimensional and multimodal. As before, median performance with mixtures is vastly improved with PLMCO-CE, and is not adversely affected with a single Gaussian. The large variances in the mean performance indicate that both PLMCO-CE and standard CE with mixtures are getting trapped in local minima. This is to be expected while using mixtures on a multimodal problem. Nevertheless, PLMCO seems to ameliorate this problem. While it is clear that performance is greatly improved in the initial stages, we cannot be as certain about final mean performance. Still, the overlap in confidence intervals is rather small, and best-case performance is vastly improved, suggesting that using PLMCO advantageously skews the distribution of algorithm performance. With a single Gaussian, we can see that we need more trials: the variance in the means is not small enough to draw valid conclusions about final performance. In all cases, significant early gains from using PLMCO are indisputable.

Figs.~\ref{fig:performanceHougen} and \ref{fig:performanceRosenbrock} compare performance on the Hougen and Rosenbrock functions. On the Hougen problem, as before, PLMCO provides significant gains in the early stages, but no significant final gains (all variants do eventually converge to the optimum). On the $n$-dimensional Rosenbrock, all algorithms perform very poorly, mainly owing to bad scaling (of curvature) of the problem. Even so, PLMCO does not seem to hurt the performance in any way, and arguably improves it slightly. While the 10-dimensional Rosenbrock seems to have been successfully used by \citet{krpo06}, careful examination reveals that the $10^5$ function evaluations were used to arrive within $0.02$ of the optimum, which is somewhat inefficient for a ten-dimensional problem.

\afterpage{\clearpage}

%

\subsection{Discussion and Conclusions}
In truth, this is not the whole picture. The relation between MCO and PL can be summarized as follows: each step of an MCO algorithm is essentially an entire PL problem. In addition, between iterations, MCO algorithms acquire more data by sampling, and this choice of samples is largely a heuristic such as adaptive importance sampling. In the PL community, the field of active learning~\citep[][]{frse97,daka05} comes close to addressing this problem, but this still does not address the iterative nature of the MCO algorithm. All of this is in addition to the three points mentioned in the text: while designing estimators for parametrized integrals, there are few or no published results in the literature \textit{vis-a-vis} tailoring such estimators for a search process. In particular, as described in the text, the authors are unaware of any analyses that consider the higher moments of random variables required to analyze extremal values. The same is the case with the consideration of moments coupling different estimators. Even when all these considerations are ignored, our results seem to \emph{strongly} imply that PL techniques improve MCO performance. That is, a bias-variance tradeoff at each iteration of the algorithm, for the individual $\theta$-by-$\theta$ estimators alone, even though it does not fully capture the process, does improve MCO performance.

In the particular case of the CE method, the use of PLMCO-CE strongly improves mean performance in the early stages of the algorithm. Moreover, these early gains are achieved without degrading final performance. This seems to bear out the main point of our PLMCO hypothesis: a straightforward application of well-known PL techniques to the individual iterations of an iterated MCO process should result in improved performance. Even if one were to examine the performance more closely, such as our investigation of median performance using a semilog plot, there are often significant gains from using PLMCO. These gains often occur when the problems are multimodal, or the data are sparse, and in general, when there is a real possibility of `overfitting' due to lack of data. 

One would expect that the phenomenon of having to choose between models and tune hyperparameters is the norm rather than the exception, that hard optimization problems will probably be multimodal, requiring the use of mixture distributions (and perhaps even more hyperparameters or model classes), and using naive MCO in such cases would indeed result in much overfitting, and consequently, erratic performance. Our version of PLMCO is extremely rudimentary: we ignore almost all the complexities of an iterative search algorithm, and yet these results seem to indicate that there are significant gains to be expected. An investigation of the more subtle nuances of such iterative search process is likely to lead to even more significant gains. The authors believe that several interesting insights await discovery.

\clearpage
\bibliographystyle{unsrtnat}
\bibliography{dgorurBib}
\end{document}